\documentclass[11pt]{article}
\usepackage{amsfonts,epsfig,latexsym,amscd}

\newcommand{\news}{\setcounter{equation}{0}\setcounter{thm}{0}}

\newcommand{\R}{{\mathbb{R}}}
\newcommand{\Z}{{\mathbb{Z}}}

\newcommand{\C}{{\mathbb{C}}}
\newcommand{\I}{{\mathbb{I}}}
\newcommand{\CP}{{\mathbb{C}}{{P}}}
\newcommand{\RP}{{\mathbb{R}{{P}}^2}}
\newcommand{\beq}{\begin{equation}}
\newcommand{\eeq}{\end{equation}}
\newcommand{\bea}{\begin{eqnarray}}
\newcommand{\eea}{\end{eqnarray}}
\newcommand{\ra}{\rightarrow}

\newcommand{\cd}{\partial}
\newcommand{\wt}{\widetilde}
\newcommand{\wh}{\widehat}
\newcommand{\ol}{\overline}
\newcommand{\sigvec}{\mbox{\boldmath{$\sigma$}}}
\newcommand{\lamvec}{\mbox{\boldmath{$\lambda$}}}
\newcommand{\tauvec}{\mbox{\boldmath{$\tau$}}}
\newcommand{\thetvec}{\mbox{\boldmath{$\theta$}}}
\newcommand{\nvec}{\mbox{\boldmath{$\wh{n}$}}}
\newcommand{\zerovec}{\mbox{\boldmath{$0$}}}
\newcommand{\xvec}{{\bf X}}
\newcommand{\xtvec}{{\wt{\bf X}}}
\newcommand{\yvec}{{\bf Y}}
\newcommand{\ytvec}{{\wt{\bf Y}}}

\newcommand{\msn}{{\rm M}_n}
\newcommand{\msnt}{\wt{\rm M}_n}
\newcommand{\msit}{\wt{\rm M}_1}
\newcommand{\msi}{{\rm M}_1}

\newcommand{\hol}{{\rm Hol}}
\newcommand{\holo}{\mathsf{Hol}}
\newcommand{\harm}{\mathsf{Harm}}
\newcommand{\Rrot}{{\cal R}}
\newcommand{\Pp}{\mathsf{P}}
\newcommand{\base}{{\cal B}^4}

\newtheorem{thm}{Theorem}
\newtheorem{lemma}[thm]{Lemma}
\newtheorem{prop}[thm]{Proposition}
\newtheorem{cor}[thm]{Corollary}

\newtheorem{remark}[thm]{Remark}
\newtheorem{conj}[thm]{Conjecture}
\newcommand{\id}{{\rm Id}}

\newcommand{\adeg}{|{\rm deg}|}

\begin{document}

\title{The $L^2$ geometry of spaces of harmonic maps $S^2\ra S^2$ and
$\RP\ra\RP$\\[30pt]}
\author{J.M. Speight\\[10pt]
{\sl Department of Pure Mathematics, University of Leeds} \\
{\sl Leeds LS2 9JT, England}\\
{\sl E-mail: j.m.speight@leeds.ac.uk}}

\date{}

\maketitle
\begin{abstract}
Harmonic maps from $S^2$ to $S^2$ are all weakly conformal, and so are
represented by rational maps. This paper presents a study of
the $L^2$ metric $\gamma$ on $\msn$, the space of degree $n$ harmonic maps 
$S^2\ra S^2$, or equivalently, the space of 
rational maps of degree $n$. It is proved that $\gamma$ is K\"ahler with 
respect 
to a certain natural complex structure on $\msn$. The case $n=1$ is
considered in detail: explicit formulae for $\gamma$ and its holomorphic
sectional, Ricci and scalar curvatures are obtained, it is shown that the
space has finite volume and diameter and codimension 2 boundary at infinity,
and a certain class of
Hamiltonian flows on $\msi$ is analyzed. It is proved that $\msnt$, the
space of absolute degree $n$ (an odd positive integer) harmonic maps
$\RP\ra\RP$, is a totally geodesic Lagrangian submanifold of $\msn$, and
that for all $n\geq 3$, $\msnt$ is geodesically incomplete. Possible
generalizations and the relevance of
these results to theoretical physics are briefly discussed.
\end{abstract}

\newpage

\section{Introduction}\news
\label{int}

In theoretical physics, one often regards harmonic maps $(M,g)\ra (N,h)$,
from a Riemannian manifold of dimension 2, as
static solutions of a so-called nonlinear $\sigma$-model on space-time
$(M\times\R,\eta)$, where $\eta=dt^2-g$ is the Lorentzian pseudometric. 
Those harmonic maps which minimize energy within their 
homotopy class are usually called ``lumps'' in this context, because 
generically their energy density is localized in lump-like structures
distributed over $M$. In many cases of interest, the homotopy classes of
maps $\phi:M\ra N$ are labelled by the topological degree of $\phi$, 
and the moduli
space of static
degree $n$ lumps, $\msn$, is a smooth, finite-dimensional manifold.
There is a natural Riemannian metric on $\msn$, namely the $L^2$ metric,
which assigns to each pair of tangent vectors $X,Y\in T_\phi\msn\subset
\Gamma(\phi^*TN)$ the inner product
\beq
\label{m1}
\gamma(X,Y)=\int_M d\mu_g\, h_\phi(X,Y),
\eeq
where $d\mu_g$ denotes the area measure on $(M,g)$. The physical 
interpretation of this metric is that it is the restriction to $\msn$ of the 
symmetric bilinear form defined by the kinetic energy functional of the 
parent
$\sigma$-model. Note that, unlike the harmonic map energy, the kinetic energy
(and hence $\gamma$) depends on $g$, not just the conformal class of $g$.

This paper presents a study of this metric in the cases $M=N=S^2$ and
$M=N=\RP$ with their canonical metrics. These cases are
convenient because one has complete, explicit parametrizations
of the harmonic maps in terms of rational functions. 
We will
focus particularly on the
simplest nontrivial case, degree 1 maps $S^2\ra S^2$, obtaining a quite 
thorough understanding of its $L^2$ geometry.
The choice $N=S^2$ or
$\RP$ is rather
natural from the stand-point of physics since the order parameters
of ferromagnets and nematic liquid crystals are $S^2$- and $\RP$-valued
respectively \cite{tre}. Previously, the algebraic topology of spaces of
rational maps has been studied by Segal \cite{seg} and Guest, Kozlowski,
Murayama and Yamaguchi \cite{gue}, and the algebraic topology of spaces
 of harmonic maps $S^2\ra S^m$ and $S^2\ra{\R}P^m$ 
by Furuta, Guest, Kotani and Ohnita \cite{fur}. The 
differential topology of spaces of harmonic maps $S^2\ra S^{2m}$ and 
$S^2\ra\CP^m$ has been studied by Bolton and Woodward \cite{bolwoo} and
Lemaire and Wood \cite{lemwoo} respectively. The present paper may be 
considered complementary to this body of work. 

The physical motivation behind this study is that $\sigma$-model lumps are in
many ways analogous to topological solitons in relativistic gauge theories,
such as BPS monopoles and abelian Higgs vortices. In the $S^2$ case, for
example, lumps attain a Bogomol'nyi type topological lower bound on energy
within their homotopy class, and consequently satisfy a first order 
``self-duality'' equation (namely, the Cauchy-Riemann equation). Manton 
conjectured \cite{man1} that the slow motion of $n$ BPS monopoles is well
approximated by geodesic flow with respect to the $L^2$ metric 
on the $n$-monopole moduli space. This conjecture was extended to lumps by
Ward \cite{war}, and has since been formulated and proved rigorously for 
monopoles and vortices by Stuart \cite{stu1,stu2}. The metric in the
case $n=2$, $M=\R^2$, $N=S^2$ was investigated numerically by
Leese \cite{lee}. So the physical 
motivation behind the present work is the hope that the $L^2$ metrics will
shed light on slow lump dynamics in the parent $\sigma$-model, as the
Atiyah-Hitchin \cite{atihit} and Samols \cite{sam} metrics have done for
monopole and vortex dynamics. Of course, they remain interesting and
natural geometric structures in their own right.

The rest of the paper is structured as follows.
Let $\msn$, $n\in\Z$, denote the space of degree $n$ harmonic maps
$S^2\ra S^2$. In section \ref{modu}, we give a simple, 
concrete proof that $(\msn,\gamma)$ is K\"ahler with respect to the complex
structure induced by a natural open inclusion $\msn\subset\CP^{2n+1}$. This
result was previously conjectured (in a rather more general setting) by
Ruback \cite{rub}, who gave a very persuasive formal argument in its favour.
In section \ref{isom} we show that the K\"ahler property, along with the
isometry group, almost completely determines the $L^2$ metric on $\msi$.
Specifically, we show that any K\"ahler metric on $\msi$ invariant under
the isometry group of $\gamma$ is determined by a single function of one
variable, rather than 21 functions of 6 variables, as for a generic metric
in 6 dimensions. An explicit formula for $\gamma$ is given, and it is
shown that $\msi$, although noncompact,
 has finite volume and diameter. It is shown also that the boundary of
$(\msi,\gamma)$ at infinity has codimension 2.

In section \ref{curv}, the curvature properties of $\msi$ are studied.
Explicit formulae for the holomorphic sectional curvatures of a certain
unitary frame and for the Ricci and scalar curvatures
 are derived. 
It is shown that the holomorphic sectional and scalar curvatures are 
unbounded above, and conjectured that the Ricci curvature is positive 
definite.
The relevance of 
these results to quantum lump dynamics is discussed. 

It is natural to regard
the K\"ahler form $\Omega$
on $\msn$ as a symplectic form and study the symplectic geometry of
$(\msn,\Omega)$. Such symplectic geometry has recently
 been used to study vortex dynamics in a 
non-relativistic version of the abelian Higgs model, for example \cite{rom}.
In section \ref{hami}, the most general physically meaningful
Hamiltonian flow on $(\msi,\Omega)$ is analyzed, and the corresponding
1 lump dynamics described.

In section \ref{rprp}, 
we address the geometry of spaces of harmonic maps $\RP\ra\RP$. 
Eells and Lemaire \cite{eellem} have shown that, if nonconstant, such maps
are classified homotopically by a certain odd positive integer, which we 
shall call the absolute degree (see section \ref{rprp} for a definition).
In section \ref{rprp}
 it is proved that $\msnt$, the space of absolute degree $n$
harmonic maps, is naturally identified with a certain totally geodesic
Lagrangian submanifold of $\msn$, where the symplectic form is again taken to
be the K\"ahler form. Further, it is shown that for all
$n\geq 3$, $\msnt$ is geodesically incomplete, while $\msit$ is compact.

Finally, in section \ref{conc} we speculate on possible generalizations of
this work. As an example, it is shown that the $L^2$ metric on the
space of degree 2 elliptic
functions is naturally K\"ahler.

\section{The K\"ahler property of $\msn$}\news
\label{modu}

By the Hopf degree theorem \cite{hop}, homotopy  classes of  continuous maps
$\phi:S^2\ra S^2$ are  labelled by their topological degree $n\in\Z$. A
well known argument of Lichnerowicz \cite{lic} (rediscovered independently
by physicists Belavin and Polyakov \cite{belpol} and Woo \cite{woo}) shows
that in the degree $n$ class the harmonic map energy satisfies
$E[\phi]\geq 2\pi|n|$, with equality if and only if $\phi$ is holomorphic
($n\geq 0$) or antiholomorphic ($n<0$). 
Since harmonic maps are by definition local extremals of $E$, 
(anti)holomorphic maps are harmonic, and furthermore, minimize energy
within their class. In fact, all harmonic maps $S^2\ra S^2$ are
(anti)holomorphic \cite{wood}.
Since degree $n$ and $-n$ maps are
trivially related by a change of orientation (on domain or codomain), we may,
and henceforth will, assume $n\geq 0$
without loss of generality. 

Introducing complex
stereographic coordinates $z,W$ on domain and codomain, the 
general degree $n$ harmonic map is
\beq
\label{m2}
W(z)=\frac{a_1+a_2z+\cdots+a_{n+1}z^n}{a_{n+2}+a_{n+3}z+\cdots+a_{2n+2}z^n}
\eeq
where $a_i\in\C$ are constants, $a_{n+1}$ and $a_{2n+2}$ do not both
vanish, and the numerator and denominator share no common roots. So $\msn$
is the space of degree $n$ rational maps. Clearly any point
$(\xi a_1,\ldots,\xi a_{2n+2})\in\C^{2n+2}$, $\xi\in\C^\times:=
\C\backslash\{0\}$,
determines the same rational map as $(a_1,\ldots,a_{2n+2})$, so one
may identify each rational map with a point in $\CP^{2n+1}$. This
gives a natural open inclusion $\msn\subset\CP^{2n+1}$ (not an identification,
since the ``no common roots'' condition removes a complex codimension 1
algebraic variety from $\CP^{2n+1}$) which we use to equip $\msn$ with a
topology and complex structure. This topology is natural in that it coincides
with the relative topology of $\msn$ in $C^0(S^2,S^2)$.
The metric of interest does {\em not}
derive from the inclusion $\msn\subset\CP^{2n+1}$, 
of course, but rather from definition (\ref{m1}).
We now establish:

\begin{thm} \label{kaehthm}
For all $n>0$, $(\msn,\gamma)$ is K\"ahler with respect to the 
complex structure induced by the open inclusion $\msn\subset\CP^{2n+1}$.
\end{thm}

\noindent
{\bf Proof:}
On the open set where $a_{2n+2}\neq 0$, we may introduce complex local 
coordinates $b^\alpha=a_\alpha/a_{2n+2}$, $\alpha=1,2,\ldots, 2n+1$. We may
always arrange that $a_{2n+2}\neq 0$ by a rotation of the codomain, so it
suffices to show that $\gamma$ is K\"ahler in this coordinate system.
Explicitly,
\beq
\label{m3}
\gamma=\gamma_{\alpha\beta}db^\alpha d\ol{b^\beta}
\eeq
where repeated indices are summed over, and
\bea
\label{m4}
\gamma_{\alpha\beta}&=&\int_\C\frac{dz d\ol{z}}{(1+|z|^2)^2}\,\,
\frac{1}{(1+|W|^2)^2}\frac{\cd W}{\cd b^\alpha}
\left(\ol{\frac{\cd W}{\cd b^\beta}}\right)
\\
\label{m5}
W&=&\frac{b_1+b_2z+\cdots+b_{n+1}z^n}{b_{n+2}+b_{n+3}z+\cdots+z^n}.
\eea
Note that $\gamma$ is manifestly Hermitian, that is, $\gamma_{\beta\alpha}
\equiv\ol{\gamma}_{{\alpha\beta}}$. Hence, we need only demonstrate that
\beq
\label{m6}
\frac{\cd\gamma_{\alpha\beta}}{\cd b^\delta}\equiv
\frac{\cd\gamma_{\delta\beta}}{\cd b^\alpha},\qquad
\frac{\cd\gamma_{\alpha\beta}}{\cd \ol{b^\delta}}\equiv
\frac{\cd\gamma_{\alpha\delta}}{\cd \ol{b^\beta}}
\eeq
for all $\alpha,\beta,\delta$ \cite{nak}. In fact (\ref{m6}) follow 
immediately from (\ref{m4}) and (\ref{m5}) provided one may interchange the
order of partial derivative and integral in $\cd\gamma_{\alpha\beta}/\cd
b^\delta$. But this is an immediate consequence of the following lemma, whose
proof is presented in the appendix:

\begin{lemma} \label{fublem}
Let $X$ be a compact Riemannian manifold, $F:X\times(-\epsilon,
\epsilon)\ra\R$ be smooth and $f:(-\epsilon,\epsilon)\ra\R$ such that
$$
f(x)=\int_X F(\cdot,x).
$$
Then
$$
f'(0)=\int_X F_2(\cdot,0)
$$
where $F_2:X\times(-\epsilon,\epsilon)\ra\R$ is the partial derivative of
$F$ with respect to the second entry.
\end{lemma}

\noindent
Applying this to the integrand of (\ref{m4}), with $x$ representing the
(real or imaginary part of) any of the coordinates $b^\alpha$, the result
is proved. $\Box$

Before specializing to the case $n=1$, we note two facts about $\msn$. First,
$(\msn,\gamma)$ is geodesically incomplete. This is a special case of a
more general result \cite{sadspe}. Second, both domain and codomain spheres
are isometric under the group of rotations and reflexions of $\R^3$,
 $O(3)\cong SO(3)\cup\ol{SO(3)}$. Here $\ol{SO(3)}$ denotes the orientation 
reversing component. The induced action of $O(3)\times O(3)$
on the set of continuous maps
$S^2\ra S^2$ decomposes $O(3)\times O(3)$ into degree preserving and degree
reversing components:
\bea
O(3)\times O(3)&\cong&
[(SO(3)\times SO(3))\cup(\ol{SO(3)}\times \ol{SO(3)})]\cup\nonumber \\
& &
[(SO(3)\times\ol{SO(3)})\cup(\ol{SO(3)}\times SO(3))].
\label{m11}
\eea
The degree preserving subgroup, call it $G$, acts isometrically on
$(\msn,\gamma)$. It is convenient to define $P:\msn\ra\msn$ such that
$P:W(z)\mapsto\ol{W(\ol{z})}$. Then $G\cong SO(3)\times SO(3)\times\Z_2$,
where $\Z_2=\{\id,P\}$. We shall denote the identity component of $G$ by
$G_0$.

\section{The metric on $\msi$}\news
\label{isom}

In the case $n=1$, 
the isometric action of $G_0\cong SO(3)\times SO(3)$ described above has 
cohomogeneity 1, that is, generic $G_0$ orbits have codimension 1. This is
most easily seen by identifying $\msi$ with $PL(2,\C)$. Note that the case
$n=1$ is special in that degree 1 rational maps are closed under
composition, so $\msi$ has a natural Lie group structure, namely that of
the M\"obius group $PL(2,\C)\cong SL(2,\C)/\Z_2$. Explicitly, one
identifies a rational map
\beq
\label{4}
\label{m12}
W:z\mapsto\frac{a_{11}z+a_{12}}{a_{21}z+a_{22}}
\eeq
with a projective equivalence class of $GL(2,\C)$ matrices,
\beq
\label{m13}
[M]=\left\{\xi\left(\begin{array}{cc}
a_{11} & a_{12} \\
a_{21} & a_{22}
\end{array}\right):\xi\in\C^\times\right\},
\eeq
noting that map composition and matrix multiplication correspond under the
identification. Then the
$PU(2)\cong SU(2)/\Z_2\cong SO(3)$ subgroup of $PL(2,\C)$ consists of
rotations of $S^2$, so in matrix language $G_0$ acts on $PL(2,\C)$ by
left and right $PU(2)$ matrix multiplication.

A particularly convenient moving coframe for $PL(2,\C)$ is defined as
follows.  Let $\tau_a,\, a=1,2,3$ be the standard Pauli matrices
\beq
\label{6}
\label{m14}
\tau_1=\left(\begin{array}{cc} 0 & 1 \\ 1 & 0 \end{array}\right),\quad
\tau_2=\left(\begin{array}{cc} 0 & -i \\ i & 0 \end{array}\right),\quad
\tau_3=\left(\begin{array}{cc} 1 & 0 \\ 0 & -1 \end{array}\right).
\eeq
Then any $[M]\in PL(2,\C)$ has a unique polar decomposition
\beq
\label{7}
\label{m15}
[M]=[U](\Lambda\I_2+\lamvec\cdot\tauvec)
\eeq
where $[U]=\{\pm U\}\in PU(2)$, $\lamvec\in\R^3$, $\lambda=|\lamvec|$,
$\Lambda=\sqrt{1+\lambda^2}$ and $\cdot$
denotes the $\R^3$ scalar product
\cite{penrin}. 
The moving coframe is $\{d\lambda_a,\sigma_a:a=1,2,3\}$, where $\sigma_a$
are the left-invariant one-forms on $PU(2)$ associated with the basis
$\{\frac{i}{2}\tau_a:a=1,2,3\}$ for $su(2)\cong T_{[\I_2]}PU(2)$. 
So $\msi\cong PU(2)\times\R^3$ as
a manifold (though {\em not} as a group). Physically, the lump parametrized 
by
$([U],\lamvec)$ should be thought of as located at 
$-\wh{\lamvec}\in S^2$ (where $\wh{\lamvec}=\lamvec/\lambda$), 
with ``sharpness'' $\lambda$
and internal orientation $[U]$. 
The action of $([L],[R])\in PU(2)\times PU(2)\cong G_0$ on $\msi$ in terms
of the polar decomposition is
\beq
\label{8}
\label{m17}
([L],[R]):([U],\lamvec)\mapsto([LUR],\Rrot\lamvec)
\eeq
where $\Rrot\in SO(3)$ is the rotation corresponding to $[R]\in PU(2)$
(explicitly, it has matrix components
 $\Rrot_{ab}=\frac{1}{2}{\rm tr}(\tau_a R^\dagger\tau_b R)$).
From this one sees that the $G_0$ action indeed has cohomogeneity 1, the
orbits being level sets of $\lambda$. The orbit space $\msi/G_0$ may be
identified with the radial curve $\Gamma=\{([\I_2],(0,0,\lambda)):\lambda\geq
0\}$ of rational maps $W_\lambda:z\mapsto\mu(\lambda)z$, where
$\mu(\lambda)=(\Lambda+\lambda)^2$. There is one exceptional orbit, namely
$\lambda=0$, which has codimension 3. 

The main aim of this section is to obtain an explicit formula for $\gamma$,
by applying the following:

\begin{prop}
\label{isomprop} Let $\tau$ be a $G$ invariant symmetric $(0,2)$ tensor
on $\msi$
which is Hermitian ($\tau(JX,JY)\equiv\tau(X,Y)$), and whose $J$-associated
2-form $\hat{\tau}$ ($\hat{\tau}(X,Y):=\tau(JX,Y)$) is closed. Then there
exists a smooth function $A:[0,\infty)\ra\R$ such that
\beq
\tau=
A_1d\lamvec\cdot d\lamvec +A_2(\lamvec\cdot d\lamvec)^2+
A_3\sigvec\cdot\sigvec
+A_4(\lamvec\cdot\sigvec)^2
+A_5\lamvec\cdot(\sigvec\times d\lamvec),
\eeq
where 
$$
A_1=A(\lambda),\quad
A_2=\frac{A(\lambda)}{1+\lambda^2}+\frac{A'(\lambda)}{\lambda},\quad
A_3=\left(\frac{1+2\lambda^2}{4}\right)A(\lambda),\quad
$$
\beq
\label{Aconstraints}
A_4=\left(\frac{1+\lambda^2}{4\lambda}\right)A'(\lambda),\quad
A_5=A(\lambda),
\eeq
$A'$ denotes the derivative of $A$, $\times$ denotes the $\R^3$ vector 
product and juxtaposition of 
covectors denotes symmetrized tensor product. 
\end{prop}

\noindent
{\bf Proof:}
We first show that the most general $G_0$ invariant symmetric $(0,2)$
tensor on $\msi$ is
\beq
\label{m18}
\tau=
A_1d\lamvec\cdot d\lamvec +A_2(\lamvec\cdot d\lamvec)^2+
A_3\sigvec\cdot\sigvec
+A_4(\lamvec\cdot\sigvec)^2
+A_5\lamvec\cdot(\sigvec\times d\lamvec)
+A_6\sigvec\cdot d\lamvec
+A_7(\lamvec\cdot d\lamvec)(\lamvec\cdot\sigvec),
\eeq
where $A_1,\ldots,A_7$ are functions of $\lambda$ only.

That such a $\tau$ is
$G_0$ invariant follows from the pulled back action of $G_0$ on our
moving coframe:
\beq
\label{9}
\label{m19}
([L],[R]):(d\lamvec,\sigvec)\mapsto (\Rrot d\lamvec,\Rrot\sigvec).
\eeq
We may prove that (\ref{m18}) is the most general $G_0$ invariant
symmetric $(0,2)$ tensor possible by means of the representation theory
of $SO(N)$. Any such tensor is uniquely determined by the 1-parameter
family of symmetric bilinear forms $\tau_\lambda:V_\lambda\oplus V_\lambda
\ra\R$, where $V_\lambda=T_{W_\lambda}\msi$, and each $\tau_\lambda$
must be invariant under the isotropy group $H_\lambda<G_0$ of $W_\lambda$.
Explicitly,
\beq
\label{m20}
H_\lambda=\left\{\begin{array}{ll}
\{([\exp(-\frac{i}{2}\psi\tau_3)],[\exp(\frac{i}{2}\psi\tau_3)]):\psi\in\R\}
\cong SO(2) & \lambda>0 \\
\{([U^\dagger],[U]):[U]\in PSU(2)\}\cong SO(3) & \lambda=0.
\end{array}\right.
\eeq
The induced action of $H_\lambda$ on $V^*_\lambda\otimes V^*_\lambda$
 leaves the
subspaces of symmetric and antisymmetric bilinear
forms invariant, that is, preserves
the splitting
\beq
\label{a3}
\label{m21}
V_\lambda^*\otimes V_\lambda^*=[V_\lambda^*\odot V_\lambda^*]\oplus
[V_\lambda^*\wedge V_\lambda^*]=:V_\lambda^{+}\oplus
V_\lambda^{-}.
\eeq
One may compute the dimension of the subspace of $V_\lambda^{+}$ on
which $H_\lambda$ acts trivially (i.e.\ the subspace of $H_\lambda$ invariant
symmetric bilinear forms) by counting the number of copies of the trivial
representation in the decomposition of $(H_\lambda,V_\lambda^{+})$ into
irreducible representations, using character orthogonality. Equation
(\ref{m18}) captures all possibilities if and only if this dimension is
7 for $\lambda>0$ and 3 for $\lambda=0$.

Consider first the generic case, $\lambda>0$, $H_\lambda\cong SO(2)$.
The $H_\lambda$ action on $V_\lambda$ has matrix representation, 
\beq
\label{a4}
R(\psi)=\left(\begin{array}{cccccc}
\cos\psi & \sin\psi & 0 & 0 & 0 & 0 \\
-\sin\psi& \cos\psi & 0 & 0 & 0 & 0 \\
0 & 0 & 1 & 0 & 0 & 0 \\
0 & 0 & 0 & \cos\psi & \sin\psi & 0 \\
0 & 0 & 0 &-\sin\psi & \cos\psi & 0 \\
0 & 0 & 0 & 0 & 0 & 1
\end{array}\right)
\eeq
relative to
the ordered basis $(\cd/\cd\lambda_1,\ldots,\theta_3)$, where $\{\theta_a\}$
are the left-invariant vector fields dual to $\{\sigma_a\}$. Hence 
the character $\chi:H_\lambda\ra\R$ of this representation is
\beq
\label{a5}
\chi(\psi)={\rm tr}\, R(\psi)=2+4\cos\psi.
\eeq
The character of the induced representation of $SO(2)$ on $V_\lambda^{\pm}$
is \cite{jon1}
\beq
\label{a6}
\wt{\chi}_\pm(\psi)=\frac{1}{2}\{[{\rm tr}\, R(\psi)]^2\pm
{\rm tr}[R(\psi)^2]\}=
\left\{\begin{array}{ll}
7+8\cos\psi+6\cos 2\psi & \mbox{symmetric} \\
5+8\cos\psi+2\cos 2\psi & \mbox{antisymmetric.}
\end{array}\right.
\eeq
We shall make use of the result for $V_\lambda^-$ when analyzing the 
$J$-associated 2-form $\hat{\tau}$.
Since $SO(N)$ is a compact Lie group, the characters of inequivalent
irreducible representations are orthogonal functions on $SO(N)$ with
respect to the Haar measure. One may therefore extract the coefficient
$a_0^\pm$ of the trivial character ($\chi_0(\psi)=1$) from the decomposition
\beq
\label{a7}
\wt{\chi}_{\pm}=\sum_n a_n^\pm\chi_n
\eeq
of $\wt{\chi}_\pm$ into irreducible representations 
by taking the character inner product of both
sides of (\ref{a7}) with $\chi_0$,
\beq
\label{a8}
a_0^\pm\int_{SO(N)}d\mu\, \chi_0^2=\int_{SO(N)}d\mu\, \chi_0\wt{\chi}_\pm
\eeq
where $d\mu$ is the Haar measure. For $SO(2)$, $d\mu=d\psi/2\pi$, so
\beq
\label{a9}
\label{13}
a_0^\pm=\int_0^{2\pi}\frac{d\psi}{2\pi}\, \wt{\chi}_\pm(\psi)=
\left\{\begin{array}{ll}
7 & \mbox{symmetric} \\
5 & \mbox{antisymmetric,}
\end{array}\right.
\eeq
in agreement with (\ref{m18}).

In the special case $\lambda=0$, the isotropy group is $H_0\cong SO(3)$
whose action on $V_\lambda$ has matrix representation
\beq
\label{a10}
R(\psi,\nvec)=\left(\begin{array}{cc}
{\cal O}(\psi,\nvec) & 0 \\
0 & {\cal O}(\psi,\nvec)
\end{array}\right)
\eeq
where $(\psi,\nvec)$ parametrizes the rotation through angle $\psi$ about
axis $\nvec\in S^2$ and ${\cal O}(\psi,\nvec)$ is the associated $SO(3)$
matrix. The character of this representation is
\beq
\label{a11}
\chi(\psi,\nvec)=2\, {\rm tr}\, {\cal O}(\psi,\nvec)=2(1+e^{i\psi}+e^{-i\psi})
=2+4\cos\psi.
\eeq
It follows from (\ref{a11}), (\ref{a5}) and (\ref{a6}) that the characters of
the induced representations on $V_\lambda^{\pm}$ are the same trigonometric
functions $\wt{\chi}_\pm(\psi)$ above, independent of $\nvec$. Once again,
we may extract $a_0^\pm$ using character orthogonality, but now we must
integrate over $SO(3)$ using the Haar measure, which is
\beq
d\mu=\frac{1}{\pi}\sin^2\frac{\psi}{2}\, d\psi
\eeq
after integrating over $\nvec$ 
\cite{jon2}. The
result is
\beq
\label{a13}
\label{16}
a_0^\pm=\frac{1}{\pi}\int_0^{2\pi}d\psi\sin^2\frac{\psi}{2}\, \wt{\chi}_\pm
(\psi)=
\left\{\begin{array}{ll}
3 & \mbox{symmetric} \\
1 & \mbox{antisymmetric,}
\end{array}\right.
\eeq
which proves the initial claim.

Since $\tau$ is $G$ invariant (not merely $G_0$ invariant), it
must also
be invariant under the discrete isometry $P$, which in matrix terms is
$P:[M]\mapsto[\ol{M}]$ (entrywise complex
conjugation). The pulled-back action on the
moving coframe is
\beq
P^*:(d\lamvec,\sigvec)\mapsto(d\lambda_1,-d\lambda_2,d\lambda_3,-\sigma_1,
\sigma_2,-\sigma_3),
\eeq
implying that $A_6\equiv A_7\equiv 0$.

It remains to show that the coefficient functions $A_1,\ldots,A_5$ are
determined by the single function $A$ as in (\ref{Aconstraints}). This 
follows from Hermiticity of $\tau$ and closure of $\hat{\tau}$.
Recall that the complex structure on $\msi$
is inherited from the open inclusion
$\msi\subset{\C}P^3$. For example, on the open set where $a_{11}\neq 0$
(equation (\ref{4})), we may use the inhomogeneous coordinates
\beq
\label{18.5}
b_1=\frac{a_{12}}{a_{11}},\quad
b_2=\frac{a_{21}}{a_{11}},\quad
b_3=\frac{a_{22}}{a_{11}},
\eeq
to define a complex coordinate chart. 
 This chart contains the curve $\Gamma$ we
are using to parametrize the orbit space $\msi/G_0$. It
is a simple matter to write down the almost complex structure $J$
associated with this complex structure, in terms of the 
basis $\{\cd/\cd\lambda_a,\theta_a:a=1,2,3\}$ for $V_\lambda$, namely
\bea
J:\frac{\cd\, }{\cd\lambda_1}\mapsto\frac{2}{\Lambda}\left(
\theta_1-\frac{\lambda}{2}\frac{\cd\,\, }{\cd\lambda_2}\right),\quad\,\,\,
J:\frac{\cd\,\, }{\cd\lambda_2}\mapsto\frac{2}{\Lambda}\left(
\theta_2+\frac{\lambda}{2}\frac{\cd\,\, }{\cd\lambda_1}\right),\quad
\,\,\,\,
J:\frac{\cd\,\, }{\cd\lambda_3}\mapsto\frac{2}{\Lambda}
\theta_3,\quad \,\,\,\,\,
&&\nonumber \\
\label{20}
J:\theta_1\mapsto-\frac{1}{2\Lambda}\left(
\frac{\cd\,\, }{\cd\lambda_1}-2\lambda\theta_2\right),\quad
J:\theta_2\mapsto-\frac{1}{2\Lambda}\left(
\frac{\cd\,\, }{\cd\lambda_2}+2\lambda\theta_1\right),\quad
J:\theta_3\mapsto-\frac{\Lambda}{2}\frac{\cd\,\, }{\cd\lambda_3}.\quad &&
\eea
We emphasize that (\ref{20}) is valid only on tangent spaces based
at points on the curve $\Gamma$. By $G$ invariance of $\tau$,
this will be all the information we need.

Hermiticity of $\tau$, $\tau_\lambda(JX,JY)\equiv\tau_\lambda(X,Y)$
for all $X,Y\in V_\lambda$,  produces two nontrivial constraints on the
coefficients $A_1,\ldots,A_5$, namely,
\beq
\label{m24}
\label{hermiticity}
A_3\equiv\frac{A_1}{4}+\frac{\lambda^2}{2}A_5,\qquad
A_1+\lambda^2A_2\equiv\frac{4}{1+\lambda^2}(A_3+\lambda^2A_4).
\eeq

Let $f\in G$, and denote by the same symbol its action on $\msi$,
$f:\msi\ra\msi$.
The 2-form $\hat{\tau}(\cdot,\cdot)=\tau(J\cdot,\cdot)$
is invariant, $f^*\hat{\tau}=\hat{\tau}$, under any holomorphic
$f\in G$ since $f^*\tau=\tau$ ($G$-invariance of $\tau$) and $df_W\circ J_W=
J_{f(W)}\circ df_W$ (holomorphicity).
Similarly, $f^*\hat{\tau}=-\hat{\tau}$ for 
antiholomorphic $f\in G$. Now each $f\in G_0$
is holomorphic, so $\hat{\tau}$ is $G_0$ invariant.
We claim that the most general $G_0$ invariant 2-form on $\msi$ is
\bea
\label{m27}
\hat{\tau}&=&
\wh{A}_1(d\lamvec\cdot\sigvec-\sigvec\cdot d\lamvec)+
\wh{A}_2((\lamvec\cdot d\lamvec)(\lamvec\cdot\sigvec)-
            (\lamvec\cdot\sigvec)(\lamvec\cdot d\lamvec))
+\nonumber \\
\label{25}
& &
\wh{A}_3\lamvec\cdot(\sigvec\times\sigvec)+
\wh{A}_4\lamvec\cdot(d\lamvec\times d\lamvec)+
\wh{A}_5(d\lamvec\cdot(\lamvec\times\sigvec)-
                    (\lamvec\times\sigvec)\cdot d\lamvec),
\eea
where $\wh{A}_1,\ldots,\wh{A}_5$ are functions of $\lambda$ only, and
juxtaposition of 1-forms indicates {\em unsymmetrized} tensor
product. Clearly,
such a 2-form is $G_0$ invariant by (\ref{m19}), and is the most general
such form possible by (\ref{a9}) and (\ref{a13}). In fact, since
$P:[M]\mapsto[\ol{M}]$ is antiholomorphic, $P^*\hat{\tau}=
-\hat{\tau}$, and we may immediately conclude that $\wh{A}_5\equiv 0$.

It is a simple matter to match $\hat{\tau}(\cdot,\cdot)$ with 
$\gamma_\lambda(J\cdot,\cdot)$ on $V_\lambda$ using (\ref{20}),
and hence
determine $\wh{A}_1,\ldots,\wh{A}_4$ in terms of $A_1,\ldots,A_5$. 
The result is
\beq
\label{m28}
\wh{A}_1=\frac{\Lambda}{2}A_1,\quad
\wh{A}_2=\frac{\Lambda}{2}A_2,\quad
\wh{A}_3=\frac{1}{4\Lambda}(A_1+4A_3),\quad
\wh{A}_4=\frac{\lambda}{\Lambda}(A_5-A_1).
\eeq
Closure of $\hat{\tau}$ 
then gives extra constraints on the metric
coefficients $A_1,\ldots,A_5$. Using the 
standard exterior differential algebra for the left-invariant 1-forms of
$SO(3)$,
\beq
d\sigma_1=\sigma_2\wedge\sigma_3,\quad
d\sigma_2=\sigma_3\wedge\sigma_1,\quad
d\sigma_3=\sigma_1\wedge\sigma_2,
\eeq
one finds that at any $W_\lambda\in\Gamma$,
\bea
d\hat{\tau}&=&
(\wh{A}_1'-\lambda\wh{A}_2)d\lambda_3\wedge(d\lambda_1\wedge\sigma_1+
                                        d\lambda_2\wedge\sigma_2)
\nonumber \\ && +
(\wh{A}_3-\wh{A}_1)(d\lambda_1\wedge\sigma_2\wedge\sigma3+
                d\lambda_2\wedge\sigma_3\wedge\sigma1+
                d\lambda_3\wedge\sigma_1\wedge\sigma2)
\nonumber \\ && +
\lambda(\wh{A}_3'-\lambda\wh{A}_2)d\lambda_3\wedge\sigma_1\wedge\sigma_2+
(\lambda\wh{A}_4'+3\wh{A}_4)d\lambda_1\wedge d\lambda_2\wedge d\lambda_3.
\eea
Hence, $d\hat{\tau}=0$ if and only if 
\beq
\label{29}
\wh{A}_1=\wh{A}_3,\quad
\wh{A}_1'=\lambda\wh{A}_2,\quad
\wh{A}_4=0,
\eeq
the last of these following from nonsingularity of 
$\hat{\tau}$ at $\lambda=0$.
Rearranging these using (\ref{m28}) and the Hermiticity constraints
(\ref{m24}), one finds that all the metric coefficients are determined
by the single smooth function $A_1=A(\lambda)$ as in (\ref{Aconstraints}).
$\Box$

\begin{cor}\label{isomcor} The $L^2$ metric on $\msi$ is
$$
\gamma=
A_1d\lamvec\cdot d\lamvec +A_2(\lamvec\cdot d\lamvec)^2+
A_3\sigvec\cdot\sigvec
+A_4(\lamvec\cdot\sigvec)^2
+A_5\lamvec\cdot(\sigvec\times d\lamvec),
$$
where $A_1,\ldots,A_5$ are functions of $\lambda$ only, determined as in
(\ref{Aconstraints}) by the single function
\beq\label{Aform}
\label{AL2def}
A=\frac{4\pi\mu[\mu^4-4\mu^2\log\mu-1]}{(\mu^2-1)^3},
\eeq
where $\mu=(\sqrt{1+\lambda^2}+\lambda)^2$. 
\end{cor}

\noindent
{\bf Proof:}
By theorem \ref{kaehthm}, $\gamma$ is Hermitian and its $J$-associated
2-form (the K\"ahler form, henceforth denoted $\Omega$, rather than
$\hat{\gamma}$) is closed. Furthermore, $\gamma$ is $G$ invariant. Hence
proposition \ref{isomprop} applies. The formula for $A$ is obtained by
computing $\gamma_\lambda(\cd/\cd\lambda_1,\cd/\cd\lambda_1)$ using 
(\ref{m1}). $\Box$

Given a tensor $\tau$ satisfying the hypotheses of proposition 
\ref{isomprop}, it is convenient to define a second coefficient  function,
$B(\lambda):=\tau_\lambda(\theta_3,\theta_3)$. Of course, $B$ is determined
by $A$, according to (\ref{Aconstraints}):
\beq
\label{Bdef}
B(\lambda)=A_3+\lambda^2A_4\equiv\frac{1+2\lambda^2}{4}A(\lambda)+
\frac{\lambda+\lambda^3}{4}A'(\lambda).
\eeq
One finds for $\tau=\gamma$, the $L^2$ metric, that
\beq
\label{BL2def}
B=\frac{4\pi\mu^2[(\mu^2+1)\log\mu-\mu^2+1]}{(\mu^2-1)^3}.
\eeq

An explicit formula for $\gamma$ has previously appeared in the physics
literature \cite{spe1}, 
although its K\"ahler property and the resulting interdependence
of the coefficient functions was not understood, nor was a rigorous 
classification of $G$ invariant tensors on $\msi$ performed. The geodesic
flow on $(\msi,\gamma)$ has been extensively studied, also in \cite{spe1},
revealing quite complicated lump dynamics. We finish this section by
examining the large $\lambda$ behaviour of $\gamma$. Specifically, we will 
prove that $(\msi,\gamma)$ has finite volume and diameter, and describe its
boundary at infinity.

\begin{thm}\label{volthm} $(\msi,\gamma)$ has
finite volume and diameter. 
\end{thm}

\noindent
{\bf Proof:}
The volume form is
\beq
{\rm vol}=\frac{\Lambda}{2}BA^2\,
d\lambda_1\wedge
d\lambda_2\wedge
d\lambda_3\wedge
\sigma_1\wedge
\sigma_2\wedge
\sigma_3
\eeq
Hence, 
\bea
{\rm Vol}(\msi,\gamma)&=&4\pi{\rm Vol}(SO(3))\int_0^\infty d\lambda\,
\lambda^2\, \frac{\sqrt{1+\lambda^2}}{2}BA^2 \nonumber \\
&=&\frac{\pi}{16}{\rm Vol}(SO(3))\int_1^\infty \frac{d\mu}{\mu}\,
\left(\mu-\frac{1}{\mu}\right)^2BA^2 \nonumber \\
&<&c+\pi^3{\rm Vol}(SO(3))\int_2^\infty d\mu\,
\mu\left(2^4\frac{\log\mu}{\mu^2}\right)\left(\frac{2^3}{\mu}\right)^2
\eea
where $c$ is a constant (the volume from $\mu=1$ to $\mu=2$). Hence
$(\msi,\gamma)$ has finite volume.

One may similarly bound the diameter of $(\msi,\gamma)$,
\beq
{\rm diam}(\msi,\gamma):=\sup_{W_1, W_2\in\msi}d(W_1,W_2).
\eeq
 By the triangle
inequality,
\beq
{\rm diam}(\msi,\gamma)\leq 2\sup_{W\in\msi}d(W,{\rm Id}).
\eeq
The
distance of any map $W$ from  ${\rm Id}$ is
bounded above by the sum of the length of the radial curve from
$([U],\lamvec)$ to $([U],\zerovec)$ and the distance in $SO(3)$ from $[U]$ to
$[\I]$ with respect to the bi-invariant metric 
$A_3(0)\sigvec\cdot\sigvec$. The
latter contribution is bounded independent of $[U]$ by compactness of $SO(3)$,
and the former is, by $G_0$ invariance, bounded above
by the length of
the curve $\Gamma$. But
\bea
{\rm length}(\Gamma)&=&\int_0^\infty d\lambda\, \sqrt{A_1+\lambda^2 A_2}=
\int_1^\infty \frac{d\mu}{\mu}\sqrt{B} \nonumber \\
&<& c+8\sqrt{\pi}\int_2^\infty d\mu \frac{\sqrt{\log\mu}}{\mu^2}<\infty.
\eea
Hence $(\msi,\gamma)$ has finite diameter. $\Box$

For both estimates, the key point is that $A(\lambda)$ and $A'(\lambda)$
decay sufficiently rapidly as $\lambda\ra\infty$ to guarantee convergence
of the integrals. Note that while every $G$
invariant K\"ahler metric on $\msi$ is determined by a single function
$A(\lambda)$, the converse is false: not every $A(\lambda)$ defines such a
metric since one must also demand that $\gamma$
 be positive definite. This places
one nontrivial constraint on $A$,
\beq
\label{39}
\gamma_\lambda(\frac{\cd\, }{\cd\lambda_3},\frac{\cd\, }{\cd\lambda_3})>0
\quad\Rightarrow\quad
\frac{A'}{A}>-\frac{1+2\lambda^2}{\lambda+\lambda^3},
\eeq
and one trivial constraint ($A>0$), which together bound the decay rate of
$A(\lambda)$ as $\lambda\ra\infty$. Integrating inequality (\ref{39})
yields, for example,
\beq
\label{40}
A(\lambda)>\frac{\sqrt{2}A(1)}{\lambda\sqrt{1+\lambda^2}}\qquad 
\forall\lambda>1,
\eeq
so the decay of $A$ cannot be faster than $O(1/\lambda\Lambda)$. It is 
interesting to note that the asymptotic behaviour of the $L^2$ metric
saturates this bound, namely $\lim_{\lambda\ra\infty}\lambda\Lambda A=\pi$.

As shown above, the boundary of $(\msi,\gamma)$ at infinity lies at finite 
proper distance, so the space is geodesically incomplete. One expects,
however, that generic geodesics do {\em not} escape to infinity, since the
boundary has codimension 2.\, In fact, this boundary may be identified with
the base space $\base$ of the fibration of generic (that is $\lambda>0$)
orbits by circles
\beq
\label{fibre}
\begin{CD}
SO(3)\times S^2 \\
@V{\pi}VV \\
\base
\end{CD}
\eeq
defined as follows: for each $\wh{\lamvec}\in S^2$, the fibre containing
$\Rrot\in SO(3)\times\{\wh{\lamvec}\}$ is the orbit of $\Rrot$ under the
isotropy group (with respect to the standard $SO(3)$ action on $S^2$) of
$\wh{\lamvec}$. To see this, identify $G_0/H_\lambda$ ($\lambda>0$) with
the $G_0$ orbit of $W_\lambda\in\Gamma$, and note that the image of the left 
invariant frame
\beq
\{(\theta_a,0),(-\theta_a,\theta_a):a=1,2,3\}
\eeq
of $G_0$ maps under the linearized coset projection at the identity to
\beq
\{\theta_1,\theta_2,\theta_3,-\lambda\frac{\cd\,\,}{\cd\lambda_2},
\lambda\frac{\cd\,\,}{\cd\lambda_1},0\}\subset V_\lambda.
\eeq
All but one of the non-zero image vectors have length bounded away from $0$
for $\lambda>0$. However, 
$||\theta_3||^2=B(\lambda)\sim\pi\log\lambda/2\lambda^4\ra 0$ as $\lambda\ra
\infty$. This is the tangent vector along the fibres of $\pi:SO(3)\times S^2
\ra\base$ defined above, so as $\lambda\ra\infty$, the fibres collapse
leaving a boundary diffeomorphic to $\base$.

\section{Curvature properties}\news
\label{curv}

\subsection{Holomorphic sectional curvature}
\label{holo}

Recall that the sectional curvature of a plane
$P\in{\rm Gr}_2(T_W\msi)$ is
\beq
\sigma(X,Y):=\langle R(X,Y)Y,X\rangle,
\eeq
where $X,Y$ are orthonormal and span $P$, $\langle\cdot,\cdot\rangle=
\gamma(\cdot,\cdot)$ 
and $R$ is the Riemann curvature tensor \cite{seccurvdef}. Recall also
 that, since
$\gamma$ is Hermitian, $\gamma(X,JX)\equiv 0$ and $||JX||\equiv ||X||$, 
so one may assign to a
 line $L\in {\rm Gr_1}(T_W\msi)$ containing $X$, $||X||=1$, the 
holomorphic sectional curvature
\beq
\hol(X):=\sigma(X,JX).
\eeq
In fact, given that $\gamma$ is K\"ahler, $\hol$ uniquely determines $\sigma$
and hence $R$ \cite{kaehol}.

We shall compute the holomorphic sectional curvature of the unitary frame
$\{e_a,Je_a:a=1,2,3\}$ for $V_\lambda$,  where
\beq
\label{43}
e_1=\frac{1}{\sqrt{A_1}}\frac{\cd\, }{\cd\lambda_1},\quad
e_2=\frac{1}{\sqrt{A_1}}\frac{\cd\, }{\cd\lambda_2},\quad
e_3=\frac{1}{\sqrt{A_1+\lambda^2A_2}}\frac{\cd\, }{\cd\lambda_3}.
\eeq
Hermiticity implies that $\hol(X)\equiv \hol(JX)$, and $G$ 
invariance implies that $\hol(e_1)\equiv\hol(e_2)$, so we shall calculate
only $\hol(e_1)$ and $\hol(e_3)$. These will vary with basepoint $W_\lambda
\in\Gamma$, and hence be functions of $\lambda$.

The simpler of the two is $\hol(e_3)$:
\bea
\hol(e_3)&=&\frac{4}{(1+\lambda^2)(A_1+\lambda^2A_2)^2}
\langle
\nabla_\frac{\cd}{\cd\lambda_3}\nabla_{\theta_3}\theta_3-
\nabla_{\theta_3}\nabla_\frac{\cd}{\cd\lambda_3}\theta_3-
\nabla_{[\frac{\cd}{\cd\lambda_3},\theta_3]}\theta_3,
\frac{\cd\, }{\cd\lambda_3}
\rangle \nonumber \\
&=&\frac{4}{(1+\lambda^2)(A_1+\lambda^2A_2)^2}\left\{
\frac{\cd\, }{\cd\lambda_3}
\langle\nabla_{\theta_3}\theta_3,\frac{\cd\, }{\cd\lambda_3}\rangle-
\langle\nabla_{\theta_3}\theta_3,\nabla_\frac{\cd}{\cd\lambda_3}
\frac{\cd\, }{\cd\lambda_3}\rangle+
||\nabla_\frac{\cd}{\cd\lambda_3}\theta_3||^2\right\} \nonumber \\
&=&
\frac{1+\lambda^2}{8B^2}\left\{\left(\frac{B'}{B}-\frac{\lambda}{1+\lambda^2}
\right)B'-B''\right\}.
\label{44}
\eea
To obtain (\ref{44}), we have used metric compatibility and torsionlessness
of $\nabla$, left $SO(3)$ invariance of $\gamma$ and the Lie algebra
$su(2)\oplus\R^3$, namely,
\beq
\left[\frac{\cd\, }{\cd\lambda_a},\frac{\cd\, }{\cd\lambda_b}\right]=
\left[\frac{\cd\, }{\cd\lambda_a},\theta_b\right]=0,\quad
[\theta_a,\theta_b]=-\epsilon_{abc}\theta_c.
\eeq
Formula (\ref{44}) may be written in terms of $A$ alone using 
(\ref{Bdef}), but the result is rather messy.

Due to the more complicated expression for $Je_1$, in comparison with
$Je_3$ (see (\ref{20})), the calculation of $\hol(e_1)$ is considerably
lengthier, though no more technically difficult. We merely record the
result, which, unlike $\hol(e_3)$, simplifies somewhat when expressed
purely in terms of $A$:
\beq
\label{46}
\hol(e_1)=\frac{1}{A^2\Lambda^2}\left\{
\frac{\lambda A+\frac{1}{2}\Lambda^2A'}{(\Lambda^2+\lambda^2)A+\lambda
\Lambda^2A'}\left(\frac{\lambda A}{\Lambda^2}+A'\right)-
\frac{2+\lambda^2}{1+\lambda^2}A-\frac{3+2\lambda^2}{2\lambda}A'\right\}.
\eeq

\vbox{
\centerline{\epsfysize=3truein
\epsfbox[50 380 570 610]{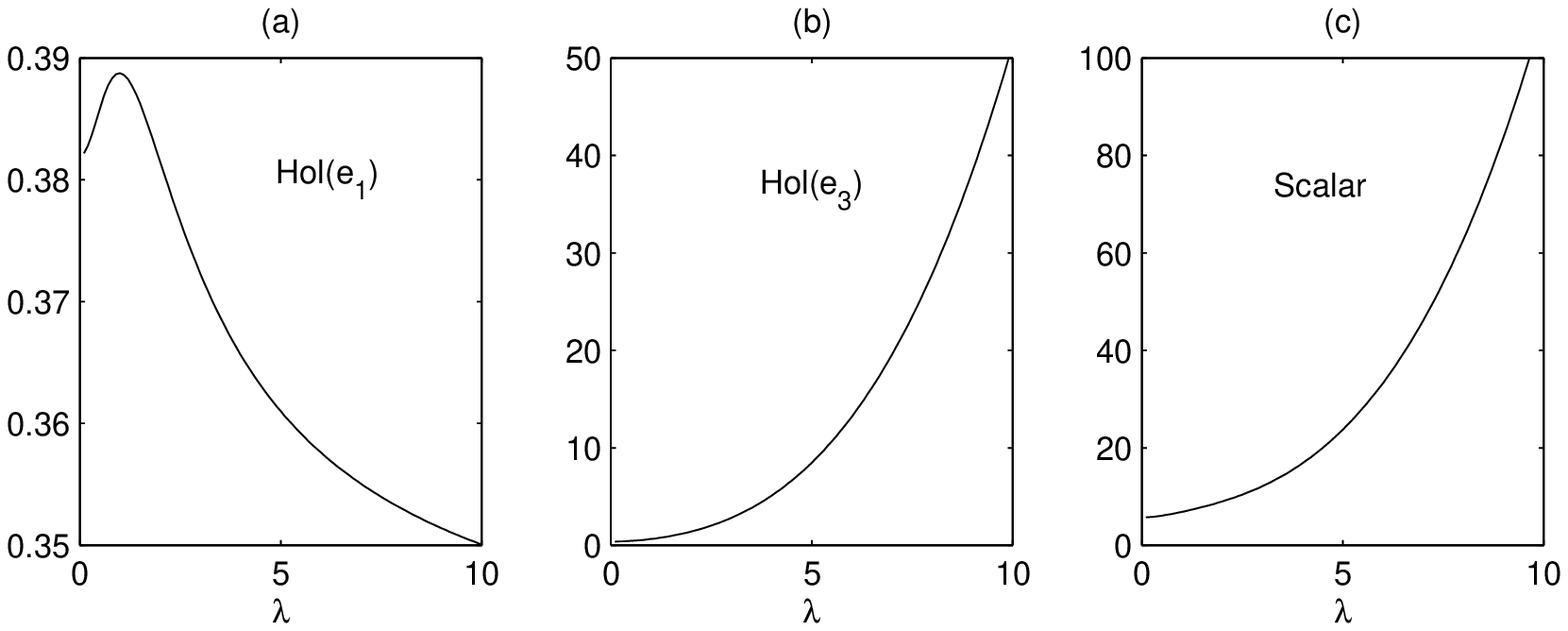}}
\noindent
{\it Figure 1: 
Plots of various curvature functions against the radial
coordinate $\lambda$ for the $L^2$ metric on $\msi$. Note the unboundedness 
of $\hol(e_3)$ and $\kappa$ (scalar curvature).}
}
\vspace{0.5cm}

Substituting the formulae (\ref{AL2def}) and (\ref{BL2def}) for $A(\lambda)$
and $B(\lambda)$ into (\ref{44}) and (\ref{46}), one obtains 
(very complicated) explicit
expressions for $\hol(e_3)$ and $\hol(e_1)$. 
Plots of these are presented in figure 1. Note that, although $\hol(e_1)$ 
is bounded, $\hol(e_3)$ is unbounded above. In fact, one finds (using
Maple, for example) that
\beq
\lim_{\lambda\ra\infty}\hol(e_1)=\frac{1}{\pi},\quad
\lim_{\lambda\ra\infty}\frac{(\log\lambda)^3}{\lambda^4}\hol(e_3)=
\frac{1}{4\pi},
\eeq
which proves:

\begin{thm}\label{holcurvthm}
The holomorphic sectional curvature of $(\msi,\gamma)$ is unbounded above.
Hence,  no isometric 
compactification of $(\msi,\gamma)$ exists, despite its finite volume and
diameter.
\end{thm}

\subsection{Ricci curvature}
\label{ricc}

Recall that the Ricci curvature $\rho$ 
of a Riemannian manifold is the symmetric $(0,2)$ tensor
\beq
\label{m50}
\rho(X,Y):={\rm tr}(V\mapsto R(V,X)Y)
\eeq
where $R$ is the Riemann curvature tensor, as before. 

\begin{prop}\label{ricciprop}
Let $\gamma$ be a $G$ invariant K\"ahler metric on $\msi$, determined
as in proposition \ref{isomprop} by the function $A$. Then the Ricci
curvature of $(\msi,\gamma)$ is
\beq
\label{m51}
\rho=
\bar{A}_1d\lamvec\cdot d\lamvec +\bar{A}_2(\lamvec\cdot d\lamvec)^2+
\bar{A}_3\sigvec\cdot\sigvec
+\bar{A}_4(\lamvec\cdot\sigvec)^2
+\bar{A}_5\lamvec\cdot(\sigvec\times d\lamvec).
\eeq
where $\bar{A}_1,\ldots,\bar{A}_5$ are functions of $\lambda$ only, 
determined as in (\ref{Aconstraints}) by the single function 
\beq
\label{m53}
\label{Abardef}
\bar{A}=
-{\frac {2\,\lambda \left (1+{\lambda }^{2}\right )\left ({ A'}\right )^{2}+
\left (9\,{\lambda }^{2}+4\right )A{ A'}+\lambda \left (1+{\lambda }^{2}
\right )
A
{ A''}+4\,A^{2}\lambda }{2\,\lambda A\left (A+2\,{\lambda }
^{2}A+\lambda { A'}+{\lambda }^{3}{ A'}\right )}}.
\eeq
\end{prop}

\noindent
{\bf Proof:} Since the $G$ action is isometric, $\rho$ is $G$ invariant.
Furthermore, since $\gamma$ is K\"ahler, $\rho(JX,JY)\equiv\rho(X,Y)$
\cite{kaeric}, and the associated Ricci form, $\hat{\rho}$
is closed \cite{kaeric2}. Hence, proposition \ref{isomprop} applies to
$\rho$ just as it applies to $\gamma$, and all the coefficient functions
are determined by $\rho_\lambda(\cd/\cd\lambda_1,\cd/\cd\lambda_1)=
\bar{A}(\lambda)$. But $\rho_\lambda(\cd/\cd\lambda_1,\cd/\cd\lambda_1)$ 
is determined by $A$ according to equation (\ref{m50}), which yields
formula (\ref{m53}). $\Box$

As with the metric, it is convenient to define the associated coefficient
function
\beq
\label{Bbardef}
\bar{B}(\lambda):=\rho_\lambda(\theta_3,\theta_3)=\bar{A}_3+\lambda^2\bar{A}_4
=
\frac{1+2\lambda^2}{4}\bar{A}(\lambda)+
\frac{\lambda+\lambda^3}{4}\bar{A}'(\lambda).
\eeq
An explicit formula for the Ricci curvature of the $L^2$ metric is
obtained by substituting (\ref{AL2def}) into (\ref{Abardef}). Unfortunately,
this formula is far too complicated to be instructive. However, it leads
us to:

\begin{conj}\label{ricciconj}
The Ricci curvature of the $L^2$ metric on $\msi$ is positive definite.
\end{conj}

\noindent
In support of this, note that, relative to the ordered basis
$(\cd/\cd\lambda_1,\theta_2,\cd/\cd\lambda_2,\theta_1,\cd/\cd\lambda_3,
\theta_3)$, the coefficient matrix of $\rho_\lambda$ is block diagonal with
blocks
\beq
\bar{A}\left(\begin{array}{cc}
1 & -\frac{\lambda}{2} \\
-\frac{\lambda}{2} & \frac{1+2\lambda^2}{4}
\end{array}\right),\quad
\bar{A}\left(\begin{array}{cc}
1 & \frac{\lambda}{2} \\
\frac{\lambda}{2} & \frac{1+2\lambda^2}{4}
\end{array}\right),\quad
\bar{B}\left(\begin{array}{cc}
\frac{4}{1+\lambda^2} & 0 \\
0 & 1
\end{array}\right),\quad
\eeq
whence it follows that $\rho_\lambda$ is positive definite if and only if
$\bar{A}(\lambda)>0$ and $\bar{B}(\lambda)>0$. Now $\bar{A}(0)=4$ and
$\bar{B}(0)=1$, so $\rho$ is certainly positive definite in a neighbourhood
of ${\rm Id}$, and
\beq
\label{asymp}
\lim_{\lambda\ra\infty}\lambda^2\bar{A}(\lambda)=4,\quad
\lim_{\lambda\ra\infty}(\log\lambda)^2\bar{B}(\lambda)=\frac{1}{8},
\eeq
so $\rho$ is asymptotically positive definite also. Convincing graphical
evidence in favour of the conjecture is presented in figure 2, which
contains plots of $\bar{A}$ and $\bar{B}$. 

\vbox{
\centerline{\epsfysize=2.5truein
\epsfbox[50 382 570 610]{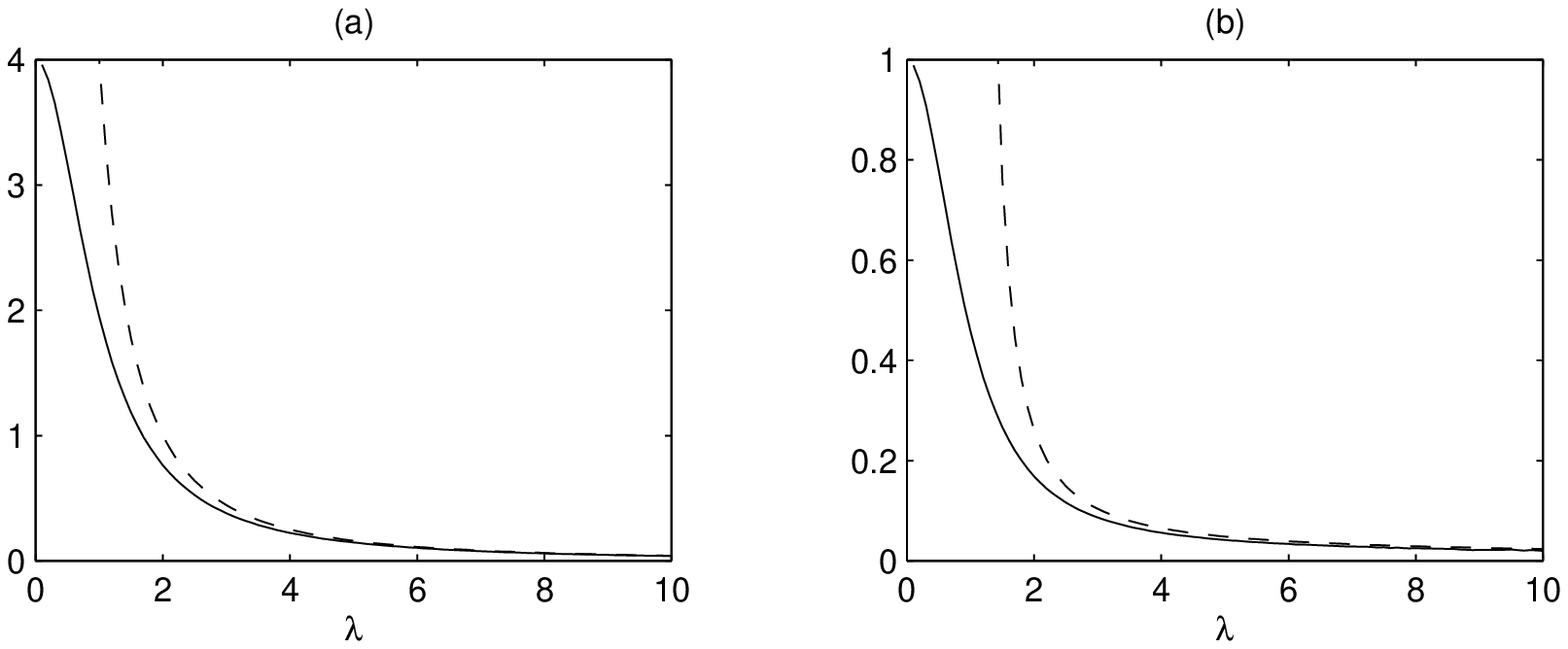}}
\noindent
{\it Figure 2: 
Plots of the coefficient functions of the Ricci curvature of $\gamma$,
(a) $\bar{A}(\lambda)$ and (b) $\bar{B}(\lambda)$. Note that both are
positive within the plot domain, and that for $\lambda\geq 4$, they are
very close to the asymptotic forms $4\lambda^{-2}$ and 
$[(\log\lambda)^{-2}]/8$ respectively (the dashed curves).}
}
\vspace{0.5cm}

We note in passing that the Einstein field equations for $G$ invariant
K\"ahler metrics,
\beq
\label{m54}
\rho=\frac{\kappa}{6}\gamma
\eeq
reduce to a single second order nonlinear ODE, explicit solutions to which
may be constructed in the Ricci flat case. The results will be described
in detail elsewhere.

\subsection{Scalar curvature}
\label{scal}

While $\hol$ and $\rho$ are
 not directly relevant to soliton dynamics, the scalar
curvature $\kappa$ certainly is, at least in the quantum regime. The standard
approach to low energy quantum $n$-soliton dynamics 
\cite{gibman} is to assume that the
quantum state is well described by a wavefunction on the $n$-soliton
moduli space $\psi:\msn\ra\C$ 
(which receives the usual probabilistic interpretation)
subject to a Schr\"odinger equation of
the form
\beq
i\frac{\cd\psi}{\cd t}=-\frac{1}{2}\Delta_\gamma\psi+V\psi
\eeq
where $\Delta_\gamma$ is the covariant Laplacian on $(\msn,\gamma)$ and
$V:\msn\ra\R$ is a potential function. The question of precisely what terms 
should be included in $V$ is somewhat controversial, and the answer likely
varies according to exact context. However, there seems to be general 
agreement that, following De Witt \cite{dew}, one should include (a positive
multiple of) $\kappa$ in $V$.
 For a recent discussion of this subject, specifically in the context of 
$\sigma$-models, see \cite{mosshi}. So the relevance of $\kappa$ to
quantum lump dynamics, as well as simple geometric curiosity, motivate us
to calculate it.

\begin{prop}\label{scalarprop}
Let $\gamma$ be a $G$ invariant K\"ahler metric on $\msi$, determined
as in proposition \ref{isomprop} by the function $A(\lambda)$. Then the
scalar curvature of $(\msi,\gamma)$ is
\beq
\label{m55}
\kappa
=4\frac{\bar{A}}{A}+2\frac{\bar{B}}{B},
\eeq
where $\bar{A}$ and $B$ are determined by $A$ as in equations (\ref{Abardef})
and (\ref{Bdef}), and
$\bar{B}$ is determined by $\bar{A}$ as in equation (\ref{Bbardef}).
\end{prop}

\noindent
{\bf Proof:}
By $G$ invariance, $\kappa$ is a function of $\lambda$ only, 
so it suffices to compute
it at $W_\lambda\in\Gamma$. Making use of the unitary frame $\{e_a,Je_a:
a=1,2,3\}$ and recalling that $\rho(JX,JY)\equiv\rho(X,Y)$, one finds
\beq
\label{kappa}
\kappa=2\sum_{a=1}^3\rho(e_a,e_a)
=2\left[2\frac{\bar{A}_1}{A}+\frac{\bar{A}_1+\lambda^2\bar{A}_2}{A_1+
\lambda_2A_2}\right],
\eeq
in the notation of proposition \ref{ricciprop}. Formula (\ref{m55}) follows
from applying the relations (\ref{hermiticity}), (\ref{Bdef}) and 
(\ref{Bbardef}) to
(\ref{kappa}). $\Box$

\begin{cor}\label{scalarcor}
The scalar curvature of the $L^2$ metric on $\msi$ is unbounded above.
\end{cor}

\noindent
{\bf Proof:} From equations (\ref{AL2def}) and (\ref{BL2def}) one has
 the limits
\beq
\label{metriclimits}
\lim_{\lambda\ra\infty}\lambda^2A(\lambda)=\pi,\quad
\lim_{\lambda\ra\infty}\frac{\lambda^4}{\log\lambda}B(\lambda)=\frac{\pi}{2},
\eeq
which together with (\ref{asymp}) and proposition \ref{scalarprop}
imply that
\beq
\label{scalarlimit}
\lim_{\lambda\ra\infty}\frac{(\log\lambda)^3}{\lambda^4}\kappa(\lambda)=
\frac{1}{2\pi}.\qquad\qquad\Box
\eeq

\begin{remark} Numerical evidence suggests that the $L^2$ metric on 
$\msi$ has strictly positive scalar curvature (see figure 1), as one would
expect, given conjecture \ref{ricciconj}.
\end{remark}

Since $(\msi,\gamma)$ is noncompact, but of finite volume, the question of
what boundary conditions to impose on the quantum wavefunction $\psi$ at
$\lambda=\infty$ when seeking bound states is non-trivial. The fact that
$\kappa\ra\infty$ as $\lambda\ra\infty$ supports the imposition of vanishing
boundary conditions for all quantum states of finite energy. One 
would expect the quantum 1-lump energy spectrum to be discrete, therefore.

\subsection{The Fubini-Study metric}
\label{fubi}

There is another natural
K\"ahler metric on $\msi$ given by the open inclusion $\msi\subset
{\C}P^3$, namely the Fubini-Study metric on ${\C}P^3$. In terms of the
local inhomogeneous coordinates 
$b_1,b_2,b_3$ (\ref{18.5}) this takes the form
\cite{hol}
\beq
\label{53}
\gamma_{FS}=\frac{(1+\sum|b_a|^2)(1+\sum db_b\ol{db}_b)-
(\sum\ol{b}_adb_a)(\sum b_b\ol{db}_b)}{(1+\sum|b_a|^2)^2}.
\eeq

\begin{prop}\label{fubstudprop}
 The Fubini-Study metric on $\msi$ is
$$
\gamma_{FS}=
A_1d\lamvec\cdot d\lamvec +A_2(\lamvec\cdot d\lamvec)^2+
A_3\sigvec\cdot\sigvec
+A_4(\lamvec\cdot\sigvec)^2
+A_5\lamvec\cdot(\sigvec\times d\lamvec),
$$
$A_1,\ldots,A_5$ being determined as in (\ref{Aconstraints})
by the single function
\beq
\label{58}
A_{FS}(\lambda)=\frac{2\mu(\lambda)}{1+\mu(\lambda)^2},
\eeq
where $\mu(\lambda)=(\sqrt{1+\lambda^2}+\lambda)^2$.
\end{prop}

\noindent
{\bf Proof:}
The isometric action of $PU(4)$ on $({\C}P^3,\gamma_{FS})$ obtained by
projecting the standard $U(4)$ action on $\C^4$ contains the $G_0$ action
on $\msi$ we have been considering. 
Furthermore, $\gamma_{FS}$ is manifestly invariant under $M\mapsto\ol{M}$
(i.e.\ $b_a\mapsto\ol{b}_a$) from (\ref{53}). Hence proposition 
\ref{isomprop} applies. It remains to compute $A_{FS}(\lambda)=
\gamma_{FS}(\cd/\cd\lambda_1,\cd/\cd\lambda_1)$ at $W_\lambda\in\Gamma$, using
(\ref{53}), which is straightforward algebra. $\Box$

Proposition \ref{fubstudprop}
 gives us several checks on our curvature calculations. It
is known that $({\C}P^3,\gamma_{FS})$ 
has constant holomorphic sectional curvature
(i.e.\ $\hol(X)$ is independent both of $X\in T_p{\C}P^3$ and base point $p$),
and is Einstein
\cite{hol}. So
substituting (\ref{58}) into (\ref{44}), (\ref{46}) and (\ref{m55}) should
yield constants. This is easily checked. One finds,
\beq
\hol_{FS}(e_1)\equiv\hol_{FS}(e_3)\equiv 4, \quad
\kappa_{FS}\equiv 48.
\eeq
Also, substituting (\ref{58}) into (\ref{m53}) demonstrates that
$\bar{A}_{FS}=8A_{FS}$, as it should.
This gives us considerable confidence in the somewhat complicated
expressions for $\hol$, $\rho$ and $\kappa$.

\section{Hamiltonian flows}\news
\label{hami}

The K\"ahler form $\Omega$ is a closed 2-form, nondegenerate
by nondegeneracy of $\gamma$, and hence a natural symplectic form on
 $\msi$. Associated with any smooth function $H:\msi\ra\R$ there is a
Hamiltonian flow, defined as the flow along the smooth vector field
$X_H$ defined such that
\beq
\label{60}
\Omega(Y,X_H)=dH(Y)
\eeq
for all vector fields $Y$. 
 Thinking
of $\msi$ as the 1-lump moduli space, only $SO(3)\times SO(3)$ invariant 
Hamiltonians make physical sense, so $H$ must be a function of $\lambda$
only. 

\begin{prop}\label{hamprop} Let $\Omega$ be the K\"ahler form associated with
a $G$ invariant K\"ahler metric on $\msi$, determined as in proposition
\ref{isomprop} by $A(\lambda)$, and $H(\lambda)$ be a
smooth, $G$ invariant function on $\msi$. 
The Hamiltonian vector field corresponding
to $(\Omega,H)$ is
\beq
X_H=\frac{2\sqrt{1+\lambda^2}H'(\lambda)}{(1+2\lambda^2)A(\lambda)
+(\lambda+\lambda^3)A'(\lambda)}
\wh{\lamvec}\cdot\thetvec.
\eeq
\end{prop}

\noindent
{\bf Proof:}
It is convenient to decompose vector fields relative to the moving
frame $\{\cd/\cd\lambda_1,\ldots,\theta_3\}$ using the notation
\beq
Y=\yvec\cdot\frac{\cd\, }{\cd\lamvec}+\ytvec\cdot\thetvec,
\eeq
that is, collecting the coefficients into a pair of $\R^3$-vector valued
functions. Recall from the proof of proposition \ref{isomprop} that the 
K\"ahler form is
\beq
\label{62}
\Omega=\wh{A}_1(d\lamvec\cdot\sigvec-\sigvec\cdot d\lamvec)+
\wh{A}_2(\lamvec\cdot d\lamvec)\wedge(\lamvec\cdot\sigvec)+
\wh{A}_1\lamvec\cdot(\sigvec\times\sigvec),
\eeq
so
 the defining equation for the Hamiltonian vector field $X_H=
\xvec\cdot\cd/\cd\lamvec+\xtvec\cdot\thetvec$ reads
\bea
\wh{A}_1(\yvec\cdot\xtvec-\ytvec\cdot\xvec)+
\wh{A}_2[(\lamvec\cdot\yvec)(\lamvec\cdot\xtvec)-
       (\lamvec\cdot\ytvec)(\lamvec\cdot\xvec)]\qquad && \nonumber \\
\qquad\qquad\qquad+
\wh{A}_1\lamvec\cdot(\ytvec\times\xtvec)&=&\frac{H'}{\lambda}\lamvec\cdot\yvec
\quad \forall\, Y \nonumber \\
\label{64}
\Rightarrow
\wh{A}_1\xvec+\wh{A}_2(\lamvec\cdot\xvec)\lamvec+\wh{A}_1\lamvec\times\xtvec
&=&\zerovec \\
\label{65}
\wh{A}_1\xtvec+\wh{A}_2(\lamvec\cdot\xtvec)\lamvec-\frac{H'}{\lambda}\lamvec
&=&\zerovec.
\eea
The pair (\ref{64}), (\ref{65}) is easily solved for $\xvec$, $\xtvec$,
yielding
\beq
X_H=\frac{H'}{\wh{A}_1+\lambda^2\wh{A}_2}\wh{\lamvec}\cdot\thetvec.
\eeq
One now uses (\ref{Aconstraints}) and (\ref{m28}) to rewrite $\wh{A}_1$ and
$\wh{A}_2$ in terms of $A$. $\Box$

Flow along $X_H$ corresponds physically to a lump which maintains
constant shape $\lambda$ and position $-\wh{\lamvec}$, while spinning 
internally at constant speed about its axis. The variation of spin speed
and sense with $\lambda$ depends on the specifics of $H(\lambda)$.

\section{The space of harmonic maps $\RP\ra\RP$}\news
\label{rprp}

We begin by recalling some relevant results of Eells and
Lemaire \cite{eellem}. The homotopy
classes of continuous maps $\phi:\RP\ra\RP$ fall into distinct families 
labelled by the induced endomorphism of the fundamental group,
$\phi_*:\pi_1(\RP)\ra\pi_1(\RP)$. Since $\pi_1(\RP)=\Z_2$, there are two
families, one for which $\phi_*$ is the zero morphism ($\phi$ maps all loops
to contractible loops), the other where $\phi_*$ is the identity morphism
($\phi$ maps noncontractible loops to noncontractible loops). The zero 
morphism family contains two classes, one of which is the trivial class. The
identity morphism family contains infinitely many classes.
Any map in this family lifts to $\wt{\phi}:S^2\ra S^2$,
\beq
\label{m80}
\begin{CD}
S^2 @>{\wt{\phi}}>> S^2 \\
@V{\pi}VV        @VV{\pi}V \\
\RP @>{\phi}>> \RP
\end{CD}
\eeq
where $\pi$ denotes the covering projection, and the different classes are
distinguished by the absolute value of the degree of $\wt{\phi}$, which may
take any {\em odd} value. We shall refer to this homotopy invariant as
the absolute degree $\adeg$ of $\phi$.

Turning to harmonic maps, all but one of the homotopy classes
described above contain harmonic representatives. Again following
\cite{eellem}, if $\phi$ belongs to the
zero morphism family, it lifts to a map $\ol{\phi}:\RP\ra S^2$ which is
also harmonic since the covering projection $\pi:S^2\ra\RP$ is a local
isometry. All harmonic maps from $\RP$ to $S^2$ are constant, so the 
nontrivial class has no harmonic representative. The moduli space of 
harmonic maps in the trivial class is thus $\RP$, and the $L^2$ metric on this
space is a constant multiple of the canonical metric. 
If $\phi$ is harmonic and belongs to the identity morphism family, it lifts
to a harmonic map $\wt{\phi}:S^2\ra S^2$ (again, because $\pi$ is a local
isometry), and the space of these is well understood
in terms of rational maps. So the task is to
identify those harmonic maps $\wt{\phi}:S^2\ra S^2$ which factor through the
quotient in (\ref{m80}). Let $p:S^2\ra S^2$ be the antipodal map
($p:z\mapsto -1/\ol{z}$ in stereographic coordinates). Then $\wt{\phi}$ 
projects to
a well defined map $\phi:\RP\ra\RP$ if and only if $\wt{\phi}\circ p=
p\circ\wt{\phi}$, or in terms of the associated rational map $W(z)$,
\beq
\label{m81}
[\ol{W(z)}]^{-1}=W(\bar{z}^{-1}).
\eeq

We now
note that given such a rational map, of degree $n>0$ say, no other degree
$n$ map projects to the same $\phi$, although $W(-1/\ol{z})$, which has 
degree $-n$, does. So we may identify $\msnt$, the moduli space of
$\adeg$ $n$ harmonic maps $\RP\ra\RP$, with the subset of $\msn$
on which (\ref{m81}) holds.

\begin{thm} $\msnt$, where $n\geq 1$ is odd,
is a totally geodesic Lagrangian submanifold of
$(\msn,\gamma,\Omega)$. 
\end{thm}

\noindent
{\bf Proof:}
Let $\Pp:\msn\ra\msn$ such that
\beq
\label{m82}
\Pp:\wt{\phi}\mapsto p\circ\wt{\phi}\circ p.
\eeq
Then $\msnt\subset\msn$ is precisely the fixed point set of $\Pp$. Since
$\Pp$ is an isometry of $(\msn,\gamma)$, in the component $\ol{SO(3)}\times
\ol{SO(3)}$, $\msnt$ is totally geodesic if it is a submanifold (i.e.\ 
nonsingular). Extending $\Pp$ naturally to ${\C}P^{2n+1}$,
one finds that
\bea
\label{m83}
&&\Pp:[a_1,\ldots,a_{n+1},a_{n+2},\ldots,a_{2n+2}]\mapsto \\
&&
[(-1)^n\ol{a}_{2n+2},(-1)^{n-1}\ol{a}_{2n+1},\ldots,-\ol{a}_{{n+3}},
\ol{a}_{{n+2}},
(-1)^{n+1}\ol{a}_{{n+1}},(-1)^n\ol{a}_{n},\ldots,\ol{a}_{2},-\ol{a}_{1}],
\nonumber
\eea
which is manifestly antiholomorphic. Hence $\Pp^*\Omega=-\Omega$, and the
K\"ahler (symplectic) form restricts to $0$ on the fixed point set. So
$\msnt$ is a Lagrangian submanifold if it is nonsingular and has (real)
dimension $2n+1$. 

It remains to check that $\msnt$ is indeed nonsingular and has half the
dimension of $\msn$. A short calculation in inhomogeneous coordinates
demonstrates that the fixed point set of $\Pp$ in ${\C}P^{2n+1}$ is smooth
with real dimension $2n+1$ if $n$ is odd, and is empty if $n$ is even (the
latter being a special case of the topological fact that no even degree
map $S^2\ra S^2$ projects to a map $\RP\ra\RP$
 in (\ref{m80})). This does not suffice for
our purposes, however, since a real codimension 2 algebraic variety must be
removed from ${\C}P^{2n+1}$ to yield $\msn$. We must verify, therefore, that
the intersection of $\msnt$ with this singular set has dimension lower than
$2n+1$.

Since the question is local, we may work in a neighbourhood of any fixed
map $\wt{\phi}$, and choose stereographic coordinates on the codomain 
which are projected from
neither $\wt{\phi}((0,0,1))$ nor $\wt{\phi}((0,0,-1))$. Then,
 in a sufficiently
small neighbourhood, all harmonic maps have rational form
\beq
\label{m84}
W(z)=\mu\frac{(z-z_1)\cdots(z-z_n)}{(z-w_1)\cdots(z-w_n)}
\eeq
where $\mu\in\C^\times$. These should be thought of as 
parametrized
by $\mu$ and a pair of {\em unordered} $n$-tuples of complex numbers
$\{w_i\}$, $\{z_i\}\in\C^n/P_n$, $P_n$ being the permutation group on
$n$ objects. Of course, in this context $\C^n/P_m\cong\C^n$ diffeomorphically
through the global coordinates $\{a_i\}$ where $(z-z_1)\cdots(z-z_n)=:
z^n+a_nz^{n-1}+\ldots+a_1$. The singular set, on which ${\rm deg}\, W<n$, is
that piece where $\{w_i\}\cap\{z_i\}\neq\emptyset$. The fixed point set of
$\Pp$ in this neighbourhood consists of maps for which
\beq
\label{m85}
\{z_1,\ldots,z_n\}=\left\{-\frac{1}{\ol{w}_1},\ldots,-\frac{1}{\ol{w}_n}
\right\}
\eeq
and
\beq
\label{m86}
|\mu|=|w_1w_2\cdots w_n|.
\eeq
Equations (\ref{m85}) and (\ref{m86}) determine a $2n+1$ dimensional
submanifold of $\C^\times\times[\C^n/P_n]\times[\C^n/P_n]$, parametrized
by $\{w_i\}\in[\C^\times]^n/P_n$ and $\arg\mu\in S^1$. From this must be 
excluded, if $n\geq 3$, the $2n-3$ dimensional variety on which $w_i=
-1/\ol{w}_j$ for some $i,j$. This still leaves a 
nonsingular $2n+1$ dimensional fixed
point set, as was to be proved. $\Box$

Note that $\Pp$ is also an antiholomorphic isometry of $\gamma_{FS},$
the Fubini-Study metric inherited from the open 
inclusion $\msn\subset{\C}P^{2n+1}$. So $\msnt$ is a totally geodesic
Lagrangian submanifold of $(\msn,\gamma_{FS},\Omega_{FS})$ also, by identical
reasoning. The metric induced on $\msnt$ by $\gamma$ is more interesting
than that induced by $\gamma_{FS}$, however, since it coincides with the
$L^2$ metric on $\msnt$. The geodesic approximation to $\RP$ lump
dynamics on $\RP$ is thus a special case of $S^2$ lump dynamics on
$S^2$. 

It is clear from the proof above that $\msnt$ is generically noncompact.
The case $n=1$ is exceptional, however. Here, as described in section
\ref{isom}, one may identify a rational map with a projective equivalence 
class $[M]$ of $GL(2,\C)$ matrices. Let $[M]$ be a fixed point of
$\Pp:PL(2,\C)\ra PL(2,\C)$. Then
\beq
\Pp:\left[\left(\begin{array}{cc}
a_{11} & a_{12} \\ a_{21} & a_{22}\end{array}\right)\right]
\mapsto
\left[\left(\begin{array}{cc}
-\ol{a}_{22} & \ol{a}_{21} \\ \ol{a}_{12} & -\ol{a}_{11}
\end{array}\right)\right]
=
\left[\left(\begin{array}{cc}
a_{11} & a_{12} \\ a_{21} & a_{22}\end{array}\right)\right].
\eeq
So there exists $\xi\in\C^\times$ such that
\beq
a_{11}=-\xi \ol{a}_{22},\quad
a_{12}= \xi \ol{a}_{21},\quad
a_{21}= \xi \ol{a}_{12},\quad
a_{22}=-\xi \ol{a}_{11},
\eeq
whence it follows that $|\xi|=1$. But then
\beq
MM^\dagger=
\left(\begin{array}{cc}
|a_{11}|^2+|a_{12}|^2 & a_{11}\ol{a}_{21}+a_{12}\ol{a}_{22} \\
a_{21}\ol{a}_{11}+a_{22}\ol{a}_{12} & |a_{21}|^2+|a_{22}|^2
\end{array}\right)
=(|a_{11}|^2+|a_{12}|^2)\I_2,
\eeq
so $[M]\in PU(2)\cong SO(3)$. Hence $\msit$ consists of the rotation orbit
of ${\rm Id}:z\mapsto z$, and the induced metric $\wt{\gamma}$ on $\msit$ is
\beq
\wt{\gamma}=A_3(0)\sigvec\cdot\sigvec,
\eeq
the standard bi-invariant metric on $SO(3)$, up to a constant factor. Each
$\phi\in\msit$ has completely uniform energy density, so it is rather
misleading to call these solutions ``$\RP$ lumps.''

For higher $n$ the possibilities are more varied. For example, the energy
density of $[z\mapsto z^n]\in\msnt$, $n\geq 3$, is concentrated in a 
symmetric band centred on a (projected) great circle on $\RP$, the band
being narrower for larger $n$. Considering rational maps of
the form (\ref{m84}), with parameters satisfying (\ref{m85}) and (\ref{m86}),
a sharp lump-like structure may be induced 
by arranging that one of the poles of $W$ be close to one of the zeroes, for
example by choosing $w_2$ close to $-1/\ol{w}_1$, while keeping the other
poles and zeroes well separated. Since lumps are associated with close
pole-zero {\em pairs}, and the poles determine the zeroes (they must be 
antipodal), for $\phi\in\msnt$ at most $(n-1)/2$ distinct lumps in
the energy distribution are possible. 

The origin of the noncompactness of $\msnt$, $n\geq 3$, is that when
$w_2\ra -1/\ol{w}_1$, say, the degree of $\wt{\phi}$ drops by 2, that is
a lump (or, in the lifted picture, an antipodal pair of lumps) forms,
collapses to an infinitely sharp spike and disappears. In fact, there are
geodesics with respect to $\wt{\gamma}$ which reach such singularities in 
finite time. We conclude this section by establishing:

\begin{thm} For all $n\geq 3$, $(\msnt,\wt{\gamma})$ is geodesically
incomplete.
\end{thm}

\noindent
{\bf Proof:} It suffices \cite{incomplete} to exhibit a curve of finite
length which converges to infinity, that is, escapes every compact subset
of $\msnt$. Consider the curve $\Gamma:[\frac{1}{2},1)\ni\rho\mapsto
W_\rho\in\msnt$ where
\beq
\label{el1}
W_\rho(z)=\rho z^{n-2}\frac{(z+1)(z-\rho^{-1})}{(z-1)(z+\rho)},
\eeq
which certainly converges to infinity (as $\rho\ra 1$). The induced metric
on $\Gamma$ is $\wt{\gamma}_\Gamma=f(\rho)d\rho^2$, where
\beq
\label{el2}
f(\rho)=\int_\C\frac{dzd\ol{z}}{(1+|z|^2)^2}\frac{1}{(1+|W_\rho|^2)^2}
\left|\frac{\cd W_\rho}{\cd\rho}\right|^2.
\eeq
We now appeal to a technical lemma whose proof is postponed to the appendix:

\begin{lemma} \label{estlem}
There exist $C>0$ and $\rho_*\in(0,1)$ such that
for all $\rho\in(\rho_*,1)$,
$$
f(\rho)<C\left[1+\log\left(\frac{1}{1-\rho}\right)\right].
$$
\end{lemma}

\noindent Hence, the length of $\Gamma$,
\beq
\int_\frac{1}{2}^1 d\rho\sqrt{f(\rho)}<C\left[1+\int_{\rho_*}^1 d\rho
\sqrt{1+\log\left(\frac{1}{1-\rho}\right)}\right]
\eeq
is finite. $\Box$

Note that this result does not follow directly from the results of
\cite{sadspe} previously mentioned (incompleteness of
$\msn$), although the method of proof is 
similar. Recall that geodesic flow on $(\msnt,\wt{\gamma})$ is conjectured to
approximate closely the low energy dynamics of the $\RP$ $\sigma$ model on
spacetime $\RP\times\R$.
So the geodesic approximation predicts that $\RP$ lumps on $\RP$ may
collapse and form singularities in finite time, just as it does for $S^2$
lumps on any compact Riemann surface. 
In fact, little is known about singularity formation in the full $(2+1)$
dimensional system, although there is some numerical evidence in favour
of lump collapse \cite{lin,piezak}.

\section{Concluding remarks}\news
\label{conc}

One could hope to generalize the results of this paper in at least two
directions. Replacing the domain 2-sphere by an arbitrary compact Riemann
surface $\Sigma$, one could study the $L^2$ metric on the space 
$\holo_n(\Sigma)$
of degree
$n$ (anti)holomorphic maps $\Sigma\ra S^2$. If nonempty, $\holo_n(\Sigma)$
is the space of {\em minimal energy} degree $n$
harmonic maps  (if
empty, for example $\holo_{\pm 1}(T^2)=\emptyset$, 
there exists no minimal energy
degree $n$ harmonic map), which is the space of most direct interest to
physicists, rather than the space of {\em all} harmonic maps. 
$\holo_n(\Sigma)$ has the structure of a complex algebraic variety, so one
would expect theorem \ref{kaehthm}, the K\"ahler property of the $L^2$
metric, to generalize to this situation. (In fact, $\holo_n(\Sigma)$ may
not be smooth if $|n|\leq 2\,
{\rm genus}(\Sigma)-2$, but the K\"ahler property
should still hold in the complement of the singular set.)

As an example, consider $\holo_2(T^2)$. It was proved in \cite{spe2} that
$\holo_2(T^2)$ is homeomorphic (in $C^0$ topology) to the complex
homogeneous space $[PL(2,\C)\times T^2]/V_4$, where $V_4$ is a certain
Viergruppe (finite group of order $4$, each element being its own inverse).
So $\holo_2(T^2)$ inherits a natural complex structure from the covering
space $PL(2,\C)\times T^2$, and it suffices to show that the lift of the
$L^2$ metric is K\"ahler. Explicitly, a point 
\beq
\left(\left[\left(\begin{array}{cc} a_{1} & a_{2} \\ a_3 & a_4
\end{array}\right)\right],s\right)\in PL(2,\C)\times T^2
\eeq
is identified with the degree 2  holomorphic map
\beq
W(z)=\frac{a_1\wp(z-s)+a_2}{a_3\wp(z-s)+a_4}
\eeq
where $\wp$ is the Weierstrass p-function. Introducing inhomogeneous
coordinates on $PL(2,\C)$, an essentially identical argument to that of the
proof of theorem \ref{kaehthm} establishes:

\begin{thm} The $L^2$ metric $\gamma$
on $\holo_2(T^2)$ is K\"ahler with respect to
the complex structure induced by the identification with
$[PL(2,\C)\times T^2]/V_4$.
\end{thm}

\noindent
It is interesting to note that $(\holo_2(T^2),\gamma)$, like $(M_1,\gamma)$
 has finite diameter, leading one to expect that theorem \ref{volthm} should
generalize to $(\holo_n(\Sigma),\gamma)$ also.

The second natural generalization would be to replace the codomain
$S^2\cong \CP^1$ by a general projective space, $\CP^N$. Lemaire and Wood
\cite{lemwoo} have shown that the space of degree $n$, energy $4\pi E$
harmonic maps $S^2\ra\CP^2$, $\harm_{n,E}(\CP^2)$ is, in
$C^j$ topology ($j\geq 2$), a disjoint union of
smooth manifolds indexed by total ramification index. Further, there is an
explicit identification between each smooth component of $\harm_{n,E}(\CP^2)$
and an appropriate space of linearly full holomorphic maps $S^2\ra\CP^2$ of
fixed degree and ramification index. So again one has a natural complex
structure on the moduli space, and again one would expect the $L^2$ metric
to be K\"ahler with respect to this structure. It is even possible that the
K\"ahler property of the $L^2$ metric
may persist when the codomain itself is {\em not}
K\"ahler. Bolton and Woodward \cite{bolwoo} have conjectured that 
$\harm_E(S^{2m})$,
the
space of energy $4\pi E$ harmonic maps $S^2\ra S^{2m}$, is a complex
algebraic variety (of dimension $2E+m^2$). If true, it would be natural to
ask again whether $\gamma$ is K\"ahler, at least on the smooth part
of $\harm_E(S^{2m})$. Note that both these generalizations lie beyond the
scope of Ruback's formal argument \cite{rub}.

\section*{Acknowledgments}

The author wishes to thank Nick Manton, Matt
Szyndel and John C. Wood for useful discussions. He holds an
EPSRC
Postdoctoral
Research Fellowship in Mathematics. 

\section*{Appendix: Proofs of Lemmas \ref{fublem} and \ref{estlem}}

\subsection*{Lemma \ref{fublem}}

Since $F:X\times (-\epsilon,\epsilon)\ra\R$ is smooth, its partial 
derivative with respect to the second entry, $F_2$ is continuous. Hence
the restriction $\wt{F}_2:X\times[0,x]\ra\R$, $0<x<\epsilon$, is integrable
(its domain is compact). Thus, by the Fubini theorem \cite{fubini},
$$
\int_X\left\{\int_{[0,x]}\wt{F}_2\right\}=
\int_{[0,x]}\left\{\int_X\wt{F}_2\right\}.
$$
But
$
\int_{[0,x]}\wt{F}_2\equiv F(\cdot,x)-F(\cdot,0),
$
so the left hand side is $f(x)-f(0)$. Hence, by the fundamental theorem of
calculus,
$$
f'(0)=\left.\left\{\int_X\wt{F}_2\right\}\right|_{x=0}=\int_X
F_2(\cdot,0).
$$

\subsection*{Lemma \ref{estlem}}

From (\ref{el1}) and (\ref{el2}) one finds that
$$
f(\rho)=\int_\C dzd\ol{z}\, F(z,\rho),
$$
where
$$
F(z,\rho)=\frac{|1+z^2|^2}{(1+|z|^2)^2}\times
\frac{|z|^{2(n-2)}|z+1|^2|z-1|^2}{(|z+\rho|^2|z-1|^2+|\rho z-1|^2|z+1|^2
|z|^{2(n-2)})^2}.
$$
Fix $\epsilon\in(0,\frac{1}{2})$, and assume that $\rho$ is close to 1, that 
is $0<\rho-1<\epsilon$. Then $F(\cdot,\rho)$ may be bounded independent of
$\rho$ except on the union of disks $D_\epsilon(-1)\cup D_\epsilon(1)$, where
one or other of the terms in the denominator may vanish (here
$D_r(z_0)=\{z\in\C:|z-z_0|<r\}$). We shall denote positive constants
(independent of $z$ and $\rho$) by $C_1$, $C_2$ etc..
On $D_\epsilon(-1)$ there exists $C_1$ such that
$$
F(z,\rho)<\frac{C_1|z+1|^2}{(|z+\rho|^2+\alpha\rho^2|z+1|^2)^2}
$$
where $\alpha=(1-\epsilon)^{2(n-2)}<1$. Hence, defining $re^{i\theta}:=
z+1$,
\begin{eqnarray*}
\int_{D_\epsilon(-1)}dzd\ol{z}\, F(z,\rho)&<&
C_1\int_0^{2\pi}d\theta\int_0^\epsilon dr\frac{r^3}{[(1+\alpha\rho^2)r^2
-2(1-\rho)\cos\theta\, r +(1-\rho)^2]^2} \\
&<& C_2\int_0^1 dr\frac{r^3}{[\frac{3}{2}r^2-2(1-\rho)r +(1-\rho)^2]^2}
\end{eqnarray*}
provided $\rho>(2\alpha)^{-1/2}\in (0,1)$. Then, rescaling $r\mapsto
r/(1-\rho)$, one finds that
\begin{eqnarray*}
\int_{D_\epsilon(-1)}dzd\ol{z}\, F(z,\rho) &<& C_2\int_0^{(1-\rho)^{-1}}
dr\frac{r^{3}}{[\frac{3}{2}r^2-2r^2+1]^2} \\
&<&C_3 +C_4\int_1^{(1-\rho)^{-1}}\frac{dr}{r} \\
&<&C_5[1-\log(1-\rho)].
\end{eqnarray*}
Noting that $\rho>(2\alpha)^{-1/2}$ implies $1+\alpha\rho^2>3\alpha/2$, 
one finds a 
similar estimate for the contribution from $D_\epsilon(1)$:
\begin{eqnarray*}
\int_{D_\epsilon(1)}dzd\ol{z}\, F(z,\rho) &<&
C_6\int_{D_\epsilon(1)}dzd\ol{z}\frac{|z-1|^2}{(|z-1|^2+\alpha|\rho z-1|^2)^2}
\\
&<& C_7\int_0^1dr\frac{r^3}{[(1+\alpha\rho^2)r^2-2\alpha\rho(1-\rho)r
+\alpha(1-\rho)^2]^2} \\
&<& C_8\int_0^{(1-\rho)^{-1}}dr\frac{r^3}{[\frac{3}{2}r^2-2r+1]^2} \\
&<& C_9[1-\log(1-\rho)]
\end{eqnarray*}
Since $F$ is bounded independent of $\rho$ on $U=\C\backslash[D_\epsilon(-1)
\cup D_\epsilon(1)]$,
$$
\int_U dzd\ol{z}\, F(z,\rho)<C_{10}+C_{11}\int_1^\infty\frac{dr}{r^{2n-1}}
<C_{12}.
$$
Defining $C=C_5+C_9+C_{12}>0$ and $\rho_*=(2\alpha)^{-1/2}\in(0,1)$, the 
lemma is proved.

\end{document}